\documentclass{amsart}
\usepackage{calc,amssymb,amsthm,amsmath,setspace}
\RequirePackage[dvipsnames,usenames]{xcolor}

\usepackage{fancyvrb}
\usepackage{hyperref}
\hypersetup{
bookmarks,
bookmarksdepth=3,
bookmarksopen,
bookmarksnumbered,
pdfstartview=FitH,
colorlinks,backref,hyperindex,
linkcolor=Sepia,
anchorcolor=BurntOrange,
citecolor=MidnightBlue,
citecolor=OliveGreen,
filecolor=BlueViolet,
menucolor=Yellow,
urlcolor=OliveGreen
}
\usepackage{xspace}
\usepackage{mabliautoref}
\usepackage{colonequals}
\input{kmacros3.sty}
\usepackage{stmaryrd}

\usepackage{verbatim}
\usepackage{enumerate}

\usepackage[normalem]{ulem}

\newcommand{\fram}{\mathfrak{m}}

\DefineVerbatimEnvironment%
  {MyVerbatim}{Verbatim}
  {formatcom=\color{Violet}}

\renewcommand{\geq}{\geqslant}
\renewcommand{\leq}{\leqslant}

\newtheorem*{caveat}{Caveat}

\begin{document}
\title{The  \emph{TestIdeals} package for \emph{Macaulay2}}
\author[Alberto F.\ Boix et al.]{Alberto F.\ Boix}
\address{Department of Mathematics, Ben-Gurion University of the Negev, Beer-Sheva 8410501, Israel}
\email{fernanal@post.bgu.ac.il}
\thanks{A.F.\,Boix was supported by Israel Science Foundation (grant No. 844/14) and Spanish Ministerio de Econom\'ia y Competitividad MTM2016-7881-P}

\author[]{Daniel J.\ Hern\'andez}
\address{Department of Mathematics, University of Kansas, Lawrence, KS~66045, USA}
\email{hernandez@ku.edu}
\thanks{D.~J.~Hern\'andez was partially supported by NSF DMS \#1600702.}

\author[]{Zhibek Kadyrsizova}
\address{School of Science and Technology, Nazarbayev University, Astana, 010000, Republic of Kazakhstan}
\email{zhibek.kadyrsizova@nu.edu.kz}
\thanks{Z. Kadyrsizova was partially supported by NSF DMS \#1401384 and the Barbour Scholarship at the University of Michigan.}

\author[]{Mordechai Katzman}
\address{Department of Pure Mathematics, University of Sheffield, Sheffield S37RH, United Kingdom}
\email{M.Katzman@sheffield.ac.uk}

\author[]{Sara Malec}
\address{Department of Mathematics, Hood College, Frederick, MD 21701}
\email{malec@hood.edu}

\author[]{Marcus Robinson}
\address{Department of Mathematics, University of Utah, Salt Lake City, UT~84112, USA}
\email{robinson@math.utah.edu}

\author[]{Karl Schwede}
\address{Department of Mathematics, University of Utah, Salt Lake City, UT~84112, USA}
\thanks{K.~Schwede was supported by NSF CAREER Grant DMS \#1252860/1501102, NSF FRG Grant DMS \#1265261/1501115, NSF grant \#1801849 and a Sloan Fellowship.}
\email{schwede@math.utah.edu}

\author[]{Daniel Smolkin}
\address{Department of Mathematics, University of Utah, Salt Lake City, UT~84112, USA}
\email{smolkin@math.utah.edu}
\thanks{D.~Smolkin was supported by NSF RTG Grant DMS \#1246989, NSF CAREER Grant DMS \#1252860/1501102, and NSF FRG Grant DMS \#1265261/1501115.}

\author[]{Pedro Teixeira}
\address{Department of Mathematics, Knox College, Galesburg, IL~61401, USA}
\email{pteixeir@knox.edu}

\author[]{Emily E.\ Witt}
\address{Department of Mathematics, University of Kansas, Lawrence, KS~66045, USA}
\email{witt@ku.edu}
\thanks{E.E.~Witt was partially supported by NSF DMS \#1623035.}
\date{\today}

\begin{abstract}
	This note describes a \emph{Macaulay2} package for computations in prime characteristic commutative algebra.  This includes Frobenius powers and roots, $p^{-e}$-linear and $p^{e}$-linear  maps,
  singularities defined in terms of these maps, different types of test ideals and modules, and ideals compatible with a given $p^{-e}$-linear map.
\end{abstract}



\maketitle

\section{Introduction}

This paper describes methods for computing objects and numerical invariants in prime characteristic commutative algebra, implemented in the package \emph{TestIdeals} for the computer algebra system \emph{Macaulay2} \cite{M2}.
A ring $R$ of prime characteristic $p$ comes equipped with the Frobenius endomorphism
\[ F: R \to R \ \text{ given by } F(x) =  x^p,\]
which is the basis for many constructions and methods.
Notably, the Frobenius endomorphism can be used to detect whether a ring is regular \cite{KunzCharacterizationsOfRegularLocalRings}, and further, to quantify how far a ring is from being regular, measuring the severity of a singularity.

In this direction, two notable applications of the Frobenius endomorphism are the theory of tight closure
(see \cite{HochsterHunekeTC1,HochsterFoundations} for an introduction)
and the resulting theory of test ideals
(see the survey \cite{SchwedeTuckerTestIdealSurvey}).  These methods are used by a wide group of commutative algebraists and algebraic geometers.

The \emph{TestIdeals} package was started during a \emph{Macaulay2} development workshop in 2012, hosted by Wake Forest University.
The package, at that time called \emph{PosChar}, aimed at providing a unified and efficient set of tools to study singularities in characteristic $p > 0$, and in particular, collecting and implementing several algorithms that had been described in research papers.
Since then, \emph{PosChar} was split into two packages, \emph{TestIdeals} and \emph{FThresholds}, and
much more functionality has been added by many contributors, during several more \emph{Macaulay2} development workshops.\footnote{Development workshops hosted by the University of California, Berkeley (2014, 2017), Boise State University (2015), and the University of Utah (2016).}

Starting at least with Kunz \cite{KunzCharacterizationsOfRegularLocalRings} and Fedder \cite{FedderFPureRat}, it has been known that the Frobenius endormorphism offers \emph{effective} methods for measuring singularities in positive characteristic.
However, the algorithms that form the basis of this package are the methods for computing Frobenius roots, $\star$-closures, and parameter test ideals.
These first appeared in \cite{KatzmanParameterTestIdealOfCMRings,BlickleMustataSmithDiscretenessAndRationalityOfFThresholds,BlickleMustataSmithFThresholdsOfHypersurfaces,KatzmanFrobeniusMapsOnInjectiveHulls}.
Another algorithm in \cite{KatzmanSchwedeAlgorithm} for computing prime ideals compatible with a given $p^{-e}$-linear map was implemented and used to produce the examples in that paper.  The methods for computing test ideals and test modules that were used implicitly in papers such as \cite{BlickleSchwedeTakagiZhang,KatzmanLyubeznikZhangOnDiscretenessAndRationality,SchwedeTuckerTestIdealFiniteMaps} became implementable via the Frobenius roots functionality.
Algorithms for computing $F$-pure thresholds from \cite{HernandezFInvariantsOfDiagonalHyp}, \cite{HernandezFPureThresholdOfBinomial}, and \cite{HernandezTeixeiraFThresholdFunctions} form some of the key methods in the forthcoming \emph{FThresholds} package.

\subsection*{Acknowledgements}
The authors thank other contributors of the \emph{TestIdeals} package, namely
Erin Bela, Juliette Bruce, Drew Ellingson, Matthew Mastroeni, and Maral Mostafazadehfard.  We also thank the referees as well as Thomas Polstra for pointing out typos in a previous version of this paper.  We also thank the referees for numerous useful comments on the package.
This paper, and the finishing touches to the \emph{TestIdeals} package, were made at the University of Utah in 2018, and the visiting authors thank the Department of Mathematics for its hospitality.  The authors also thank the referees for their numerous useful comments on both this paper and the code.  Additional work on the package was also done at the Institute of Mathematics and its Applications for its generous support for a 2019 Coding Sprint on the \emph{FThresholds} package.

\section{Frobenius powers and Frobenius roots}\label{Section: Frobenius powers and Frobenius roots}

Let $R$ denote a commutative ring of prime characteristic $p$.

\begin{definition}
Given an ideal $I\subseteq R$ and an integer $e\geq 0$, we define the \emph{$p^e$-th Frobenius power of $I$}, denoted $I^{[p^e]}$, to be the ideal
generated by the $p^e$-th powers of all elements of $I$.
\end{definition}

It is easy to see that if $I$ is generated by $g_1, \dots, g_\ell$, then $I^{[p^e]}$ is generated by $g_1^{p^e}, \dots, g_\ell^{p^e}$.

\begin{definition}
Given an ideal $I\subseteq R$ and an integer $e\geq 0$, we define the \emph{$p^e$-th Frobenius root of $I$}, denoted $I^{[1/p^{e}]}$, to be the smallest ideal $J$ such that $I\subseteq J^{[p^e]}$, if such an ideal exists.
\end{definition}

Frobenius roots always exist in polynomial and power series rings, and in $F$-finite regular rings
(cf.~\cite[\S 2]{BlickleMustataSmithDiscretenessAndRationalityOfFThresholds} and \cite[\S 5]{KatzmanParameterTestIdealOfCMRings}).
Below is an example of a computation of a Frobenius root in a polynomial ring.
In \autoref{ss: math behind} we describe the main ideas behind such computations.

\medskip
{\small
\setstretch{.67}
\begin{MyVerbatim}
i1 : R = ZZ/5[x,y,z];

i2 : I = ideal(x^6*y*z + x^2*y^12*z^3 + x*y*z^18);

o2 : Ideal of R

i3 : frobeniusPower(1/5, I)

                2   3
o3 = ideal (x, y , z )

o3 : Ideal of R
\end{MyVerbatim}
}
\medskip

\subsection{The mathematics behind the computation of Frobenius roots}
\label{ss: math behind}

We can also describe Frobenius roots as follows:  In a (sufficiently local) regular ring, we have an identification of $R$ with its canonical module $\omega_R$.
On the other hand, the Grothendieck dual to the $e$-iterated Frobenius map\footnote{In a reduced ring $R$, it is often convenient to express the Frobenius endomorphism $R \to R$ as the inclusion $R \hookrightarrow R^{1/p}$, so that the source and target  are distinguished.} $R \to R^{1/p^e}$ provides a map
\begin{equation}
\label{eq.DualToFrobenius}
T : \omega_{R^{1/p^e}} \to \omega_R,
\end{equation}
which using the identification $R \cong \omega_R$ gives us a map
\[
T : R^{1/p^e} \to R.
\]
It is not difficult to prove that
$T(I^{1/p^e})=I^{[1/p^e]}$,
where $I^{1/p^e} \subseteq R^{1/p^e}$ is simply the ideal of $p^e$-th roots of the elements of $I$.
In the case that
\[
R = \mathbb{K}[x_1, \ldots, x_d],
\]
where $\mathbb{K}$ is a perfect field,
$R^{1/p^e}$ is a free $R$-module with basis consisting of the monomials
\[
\mathbf{x}^{\boldsymbol{\lambda}/p^e}= x_1^{\lambda_1/p^e} \cdots x_d^{\lambda_d/p^e},
\]
where $0 \leq \lambda_i \leq p^e-1$.
Furthermore, the map $T : R^{1/p^e} \to R$ is simply the projection defined as follows:
\[
T(x_1^{\lambda_1/p^e} \cdots x_d^{\lambda_d/p^e}) = \left\{ \begin{array}{rl} 1 &\text{if } \lambda_i = p^{e}-1 \text{ for all }i \\ 0 & \text{otherwise}  \end{array} \right.
\]
Using this, it is not difficult to see that if
\begin{equation}
\label{eqn1}
f^{1/p^e} = \sum_{\boldsymbol{\lambda}} f_{\boldsymbol{\lambda}} \mathbf{x}^{\boldsymbol{\lambda} / p^e},
\end{equation}
where $\boldsymbol{\lambda}$ runs over the tuples $(\lambda_1, \dots, \lambda_d)$ with $0 \leq \lambda_j \leq p^e-1$, and $f_{\boldsymbol{\lambda}} \in R$, then
$(f)^{[1/p^e]} = T( (f)^{1/p^e})$ is the ideal generated by the $f_{\boldsymbol{\lambda}}$.
We can then compute $I^{[1/p^e]}$ for more general $I$ by linearity: if $I=(f_1,\ldots,f_m)$, then $I^{[1/p^e]}= \sum(f_i)^{[1/p^e]}$.

\subsubsection{Complexity}
The computation of Frobenius roots is the workhorse behind many of the methods in \emph{TestIdeals}.
Hence it is important to understand how this is implemented, and its computational complexity.

The computation of Frobenius roots of ideals reduces to the case of principal ideals, and its complexity grows linearly with the number of generators of the ideal.
Furthermore the calculation of  $(f)^{[1/p^e]}$ reduces to finding the terms in the right hand side of (\ref{eqn1}), which essentially amounts to
taking each term in $f$, computing the $p^e$-th root of its coefficient, and dividing the monomial exponent vector by $p^e$ with remainder.
Hence the complexity of computing $(f)^{[1/p^e]}$ is proportional to the number of terms in $f$, and is independent of its degree.
We emphasize that the calculation of Frobenius roots \emph{does not involve the calculation of Gr\"obner bases.}

\subsection{Dual to Frobenius on quotient rings and Frobenius on top local cohomology}
\label{subsec.DualToFrobeniusOnQuotientRings}
For \emph{any} reduced ring $R$ with finite $e$-iterated Frobenius map identified with $R \hookrightarrow R^{1/p^e}$, we always have a map
\[
T_R : \omega_{R^{1/p^e}} \to \omega_R.
\]
At a maximal ideal $\fram$ of height $d$, this map is Matlis dual to the Frobenius map on $H^d_{\fram}(R)$.
Hence we can study local cohomology by studying this map.
If $R = S/I$, where $S$ is a polynomial ring over a finite field, we can implement this map \cite{FedderFPureRat, KatzmanParameterTestIdealOfCMRings}.  We briefly explain the case of a hypersurface here.

If $I = (f)$, then again $\omega_R \cong R = S/I$, and the action of $T_R$ can be computed on $S$.  In fact, if $\overline{J} \subseteq R \cong \omega_R$ has preimage  $J \subseteq S$,
then setting $u = f^{p^e-1}$, we have that $T_R({\overline J}^{1/p^e})$ is the image of $(u J)^{[1/p^e]}$ in $R$.
More generally, for non-hypersurfaces, the analog to $u$ is chosen to be an element that maps onto the generator of the cyclic $R$-module
\[
\frac{(I^{[p^e]} : I) \cap (\Omega^{[p^e]} : \Omega)}{ I^{[p^e]}},
\]
where $\Omega$ is an ideal of $S$ with $I\subseteq \Omega$ and whose image in $R$ is isomorphic to the canonical module of $R$; see \cite{KatzmanParameterTestIdealOfCMRings} for details.
Here we compute two examples.

\medskip
{\small\setstretch{.67}
\begin{MyVerbatim}
i1 : S = ZZ/5[x,y,z];

i2 : f = x^3 + y^3 + z^3;

i3 : u = f^(5 - 1);

i4 : frobeniusRoot(1, ideal u)

o4 = ideal (z, y, x)

o4 : Ideal of S

i5 : S = ZZ/7[x,y,z];

i6 : f = x^3 + y^3 + z^3;

i7 : u = f^(7 - 1);

i8 : frobeniusRoot(1, ideal u)

o8 = ideal 1

o8 : Ideal of S
\end{MyVerbatim}
}\medskip

The above example has shown that the Frobenius map on top local cohomology of the cone over the Fermat elliptic curve is injective in characteristic $7$ (the dual map $\omega_{R^{1/p}} \to \omega_R$ is surjective) and not injective in characteristic $5$ (the dual map $\omega_{R^{1/p}} \to \omega_R$ is not surjective).

\subsection{A generalization of Frobenius powers and roots}

Let $R$ be an $F$-finite regular ring.
We can extend the definition of Frobenius powers as follows.

\begin{definition}[\cite{HernandezTeixeiraWittFrobeniusPowers}]
Let  $I\subseteq R$ be an ideal.
\begin{enumerate}
 \item[(a)] If $n$ is a positive integer with base $p$ expansion  $n=d_0 + d_1 p +  \dots + d_r p^r$, we define
\[ I^{[n]}=I^{d_0} (I^{d_1})^{[p]} \cdots  (I^{d_r})^{[p^r]}.\]
 \item[(b)] If $t$ is a nonnegative rational number of the form $t = a/p^e$, we define  $I^{[t]} = {\big(I^{[a]}\big)}^{[1/p^e]}.$
 \item[(c)] If $t$ is any nonnegative rational number, and $\{a_n/p^{e_n}\}_{n\geq 1}$ is a sequence of rational numbers converging to $t$ from above, we define $I^{[t]}$
 to be the stable value of the non-decreasing chain of ideals $\{I^{[a_n/p^{e_n}]}\}_{n\geq 1}$.
\end{enumerate}
\end{definition}

\medskip
{\small\setstretch{.67}
\begin{MyVerbatim}
i1 : p = 3;

i2 : R = ZZ/p[x,y];

i3 : m = monomialIdeal(x, y);

o3 : MonomialIdeal of R

i4 : I = m^5;

o4 : MonomialIdeal of R

i5 : t = 3/5 - 1/(5*p^3);

i6 : frobeniusPower(t, I)

             2        2
o6 = ideal (y , x*y, x )

o6 : Ideal of R

i7 : frobeniusPower(t - 1/p^5, I)

o7 = ideal (y, x)

o7 : Ideal of R
\end{MyVerbatim}
}
\medskip

The generalized Frobenius powers of an ideal are test ideals (see \autoref{Section: Test Ideals}) of sufficiently general linear combinations of generators of that ideal.
Thus, the ideals computed above are test ideals of generic quintics in two variables in characteristic $p=3$, and these computations suggest that $\frac35-\frac1{5p^3}$ is a higher $F$-jumping exponent of such quintics, which is indeed the case for any $p\equiv 3 \bmod 5$ (see \cite{hernandez+etal.frobenius_examples}).
To illustrate, compare the computations above with the following.

\medskip
{\small\setstretch{.67}
\begin{MyVerbatim}
i8 : S = ZZ/3[a..f,x,y];

i9 : G = a*x^5 + b*x^4*y + c*x^3*y^2 + d*x^2*y^3 + e*x*y^4 + f*y^5;

i10 : testIdeal(t, G)

              2        2
o10 = ideal (y , x*y, x )

o10 : Ideal of S

i11 : testIdeal(t - 1/p^5, G)

o11 = ideal (y, x)

o11 : Ideal of S
\end{MyVerbatim}
}
\medskip

\begin{caveat}
   The computations of Frobenius powers with exponents whose denominators are not powers of the characteristic $p$ are often very slow, and in several instances it is more efficient to introduce more variables and compute the test ideal of a ``very general'' linear combination of the generators of the ideal, as above.
   The code used for computations of Frobenius powers needs to undergo significant optimization.
\end{caveat}

\subsection{Frobenius  roots of submodules of free modules}

Let $S$ be a polynomial ring or power series ring.
Given a submodule $M$ of the free module $S^k$,
%
%
there is a smallest submodule $N$ of $S^k$ that contains $M$, for which $M\subseteq N^{[p^e]}$.
Here,  $N^{[p^e]}$ is the submodule  of $S^k$ generated by the vectors of $N$ with all coordinates raised to the $p^e$-th power. (Cf.~\cite{KatzmanZhangAlgorithm}).
Extending the previous definitions,
we call this $N$ the $p^e$-th Frobenius root of $M$, and denote it by $M^{[1/p^e]}$.

\begin{example}
Let $R=\ZZ/p\ZZ[a,b,c,d]$, and consider the ideals $\mathfrak{m}=(a,b,c,d)$ and
$$I= (a,b) \cap (a,c) \cap (c,d) \cap (c+d, a^3+b d^2).$$
Then $R/I$ is a $3$-dimensional non-Cohen--Macaulay ring.
Matlis Duality applied to $H^2_{\mathfrak{m}} (R/I)$ with its natural Frobenius map
yields a $p^{-1}$-linear map $U$ on $\Ext^2_R(R/ I, R)$.

\medskip
{\small\setstretch{.67}
\begin{MyVerbatim}
i1 : R = ZZ/2[a,b,c,d];

i2 : I = intersect(ideal(a,b),ideal(a,c),ideal(c,d),ideal(c+d,a^3+b*d^2));

o2 : Ideal of R

i3 : f = inducedMap(R^1/I, R^1/frobenius(I));

o3 : Matrix

i4 : E2 = Ext^2(f, R^1)

o4 = {-8}  | a4+abc2+abcd a2b            |
     {-10} | a2cd3        a3cd+a3d2+bcd3 |

o4 : Matrix

i5 : target E2

o5 = cokernel {-8}  | 0     a2 b2c2 |
              {-10} | c2+d2 d4 a4d2 |

                            2
o5 : R-module, quotient of R

i6 : source E2

o6 = cokernel {-4} | 0   a  bc  |
              {-5} | c+d d2 a2d |

                            2
o6 : R-module, quotient of R
\end{MyVerbatim}
}
\medskip

The Frobenius map on $H^2_{\mathfrak{m}} (R/I)$
is injective if and only the image of
$(\Image U)^{[1/p]}$ in $\Ext^2_R(R/ I, R)$
as computed above is the whole of $R^2$.

\medskip
{\small\setstretch{.67}
\begin{MyVerbatim}
i7 : U = matrix entries E2;

             2       2
o7 : Matrix R  <--- R

i8 : A = image matrix entries relations source E2;

i9 : frobeniusRoot(1, image U)

o9 = image {-2} | 1 0 0 |
           {0}  | 0 d a |

                             2
o9 : R-module, submodule of R
\end{MyVerbatim}
}\medskip

The calculation above shows that $R/I$ is not $F$-injective (cf. \autoref{Section: F-singularities}).
Moreover, it shows that $H^2_{\mathfrak{m}} (R/I)\cong \Ann_{E^2} A^t$,
where $E$ is the injective hull of $R_{\mathfrak{m}}/ {\mathfrak{m}} R_{\mathfrak{m}}$,
and the Frobenius map on  $H^2_{\mathfrak{m}} (R/I)$ induced from Frobenius on $R$
is given by $U^t \Theta$, where $\Theta$ is the induced Frobenius on $E$.
The submodule of nilpotent elements in  $H^2_{\mathfrak{m}} (R/I)$ is
given by $\Ann_{E^2} B^t$, where $B$ is the smallest submodule of $R^2$ containing $\Image A + \Image U$ such that $U  B \subseteq B^{[p]}$.
The method \texttt{ascendModule} can be used to calculate  $B$ (see a detailed description of the similar method \texttt{ascendIdeal} in \autoref{Section: Test Ideals}).

\medskip
{\small\setstretch{.67}
\begin{MyVerbatim}
i10 : B = ascendModule(1, A, U)

o10 = image | 0   a  bc  |
            | c+d d2 a2d |

                              2
o10 : R-module, submodule of R
\end{MyVerbatim}
}\medskip

\end{example}

\section{$F$-singularities}\label{Section: F-singularities}

The \emph{TestIdeals} package includes methods for determining if a ring is $F$-injective, $F$-pure, $F$-rational, or $F$-regular.

\subsection{$F$-injectivity}

\begin{definition}
A local ring $(R, \mathfrak{m})$ is called \emph{F-injective} if the map
$H^{i}_{\mathfrak{m}}(R) \rightarrow H^{i}_{\mathfrak{m}}(R^{1/p})$ is
injective for all $i >0$. An arbitrary ring is called \emph{$F$-injective} if its
localization at each prime ideal of $R$ is $F$-injective.
\end{definition}

The function \texttt{isFInjective} determines whether the ring $R = S/I$ is
$F$-injective, where $S$ is a polynomial ring.

\medskip
{\small\setstretch{.67}
\begin{MyVerbatim}
i1 : R = ZZ/7[x,y,z]/(x^3 + y^3 + z^3);

i2 : isFInjective R

o2 = true

i3 : R = ZZ/5[x,y,z]/(x^3 + y^3 + z^3);

i4 : isFInjective R

o4 = false
\end{MyVerbatim}
}\medskip

Equivalently, a ring is $F$-injective if the maps on the cohomology of the dualizing complex
\[
h^{-i} \omega_{R^{1/p}}^{\mydot} \to h^{-i} \omega_R^{\mydot}
\]
surject for all $i$.
Note that $h^{-i} \omega_R^{\mydot} \cong \Ext^{\dim S - i}(R, S)$, the latter of which \emph{Macaulay2} readily computes.
The algorithm used by \texttt{isFInjective} works by checking the surjectivity of the dual Frobenius map
\[
\Ext^{\dim S - i}(R^{1/p}, S) \to \Ext^{\dim S - i}(R, S).
\]
We begin by computing the map $R
\rightarrow R^{1/p}$ using the \emph{PushForward} package \cite{PushForward}.
Next the algorithm computes
$\Ext^{i}( \blank, S)$ applied to the map from the previous step.  Then $R$ is $F$-injective precisely when the
cokernel of $\Ext^{i}( \blank, S)$ is trivial for $i$.

The Frobenius action on top local cohomology (dual to $\omega_R$) is usually computed in a different (faster) way than the other cohomologies, and this is modified by the \texttt{CanonicalStrategy} option.
The default value for this option is \texttt{Katzman}, which instead of using the \emph{PushForward} package, relies on the fact that we already know how to compute the Frobenius action on the canonical module, as described in \autoref{subsec.DualToFrobeniusOnQuotientRings}.

The performance of the algorithm can be improved if the ring of interest is
nice enough. If the ring is Cohen--Macaulay, then setting \texttt{AssumeCM =>
true} (the default value is \texttt{false}) lets the algorithm check the Frobenius action only on top cohomology
(which is typically much faster, as explained above).
When studying a reduced ring,  setting \texttt{AssumeReduced => true} (the default value) avoids
computing the bottom local cohomology, and when studying a normal ring, setting
\texttt{AssumeNormal => true} (the default is \texttt{false}) avoids computing the bottom two local
cohomologies.

By default the algorithm checks for $F$-injectivity at all points of $\Spec R$.  However, one
can choose to check $F$-injectivity only at the origin by setting the
option \texttt{AtOrigin} to  \texttt{true}.

\medskip
{\small\setstretch{.67}
\begin{MyVerbatim}
i1 : R = ZZ/7[x,y,z]/((x - 1)^5 + (y + 1)^5 + z^5);

i2 : isFInjective R -- R is not globally F-injective...

o2 = false

i3 : isFInjective(R, AtOrigin => true) -- but is F-injective at the origin

o3 = true
\end{MyVerbatim}
}\medskip

\subsection{$F$-regularity}

\begin{definition}[\cite{HochsterHunekeTC1,HaraWatanabeFRegFPure}]
A ring $R$ is called \emph{strongly $F$-regular} if the (big) test ideal $\tau(R)$ is equal to $R$.  Likewise a pair $(R, f^t)$ is called \emph{strongly $F$-regular} if $\tau(R, f^t) = R$.
\end{definition}

The command \texttt{isFRegular} checks whether a ring or pair is strongly
$F$-regular. Below are two examples, one of which is $F$-regular and the other one is not.

\medskip
{\small\setstretch{.67}
\begin{MyVerbatim}
i1 : R = ZZ/5[x,y,z]/(x^2 + y*z);

i2 : isFRegular R

o2 = true

i3 : R = ZZ/7[x,y,z]/(x^3 + y^3 + z^3);

i4 : isFRegular R

o4 = false
\end{MyVerbatim}
}\medskip

We can also check whether a pair $(R, f^t)$ is $F$-regular.

\medskip
{\small\setstretch{.67}
\begin{MyVerbatim}
i1 : R = ZZ/5[x,y];

i2 : f = y^2 - x^3;

i3 : isFRegular(1/2, f)

o3 = true

i4 : isFRegular(5/6, f)

o4 = false

i5 : isFRegular(4/5, f)

o5 = false

i6 : isFRegular(4/5 - 1/100000, f)

o6 = true
\end{MyVerbatim}
}\medskip

All of these checks are done by actually computing the test ideal, as described in \autoref{Section: Test Ideals}.

If the input ring is $\mathbb{Q}$-Gorenstein, then in each of the cases above, the output is a boolean indicating if the
ring is strongly $F$-regular. If the input ring is not
$\mathbb{Q}$-Gorenstein, then the algorithm can be used to determine if a
ring is strongly $F$-regular, but cannot prove that a ring is not strongly
$F$-regular (this latter functionality can, however, be enabled by setting \texttt{QGorensteinIndex => infinity}).

In the case that $R$ is $\mathbb{Q}$-Gorenstein, the algorithm works by
computing the test ideal $\tau$ of the ring  (or the pair) using \texttt{testIdeal}
and checking whether $\tau=R$.
In the non-$\mathbb{Q}$-Gorenstein case, the algorithm checks for
strong $F$-regularity by computing better and
better approximations of the test ideal, and checking whether any of these is the unit ideal.
To compute approximations of the
test ideal, the algorithm computes a test element $c$ with \texttt{testElement}
and then uses \texttt{frobeniusRoot} to compute the $e$-th root of
$c(I^{[p^{e}]} : I)$; appropriate modifications are made for pairs. If at
any step the approximation is the unit ideal, then then the algorithm
returns \texttt{true}. Otherwise the algorithm continues checking for each $e$
until a specified limit is reached. The default limit is 2, and can be
changed using the option \texttt{DepthOfSearch}.

A number of options can be used to speed up the performance of some of the
internal functions. The option \texttt{AssumeDomain} can be used if $R$ is an
integral domain, \texttt{FrobeniusRootStrategy} chooses a strategy for
internal \texttt{frobeniusRoot} calls, \texttt{MaxCartierIndex} sets the
maximum Gorenstein index to search for when working with a
$\mathbb{Q}$-Gorenstein ambient ring, and \texttt{QGorensteinIndex}
allows the user to specify the $\mathbb{Q}$-Gorenstein index of the ring.

The default behavior of \texttt{isFRegular} is that it checks for strong $F$-regularity
globally. If the option \texttt{AtOrigin} is set to \texttt{true}, the algorithm will only
check $F$-regularity at the origin, by checking whether the computed test ideal is in the irrelevant ideal.
Below are examples for both a ring and a pair.

\medskip
{\small\setstretch{.67}
\begin{MyVerbatim}
i1 : R = ZZ/7[x,y,z]/((x - 1)^3 + (y + 1)^3 + z^3);

i2 : isFRegular R -- R is not globally F-regular...

o2 = false

i3 : isFRegular(R, AtOrigin => true) -- but is F-regular at the origin

o3 = true

i4 : R = ZZ/13[x,y];

i5 : f = (y - 2)^2 - (x - 3)^3;

i6 : isFRegular(5/6, f) -- (R,f^(5/6)) is not F-regular...

o6 = false

i7 : isFRegular(5/6, f, AtOrigin => true) -- but is F-regular at the origin

o7 = true
\end{MyVerbatim}
}\medskip

\subsection{$F$-purity}
\begin{definition}
A ring $R$ is called \emph{$F$-pure} if
the inclusion $R \hookrightarrow R^{1/p^{e}}$ is a pure map, i.e.,
the tensor of this map with any $R$-module remains injective.  If $R^{1/p}$ is a finite $R$-module, this is equivalent to requiring that the inclusion $R \hookrightarrow R^{1/p}$ split as a map of $R$-modules.
\end{definition}

The function \texttt{isFPure} checks whether a ring is $F$-pure.
Either a ring or a defining ideal can be input, as seen in the following example.

\medskip
{\small\setstretch{.67}
\begin{MyVerbatim}
i1 : R = ZZ/5[x,y,z]/(x^2 + y*z);

i2 : isFPure R

o2 = true

i3 : R = ZZ/7[x,y,z]/(x^3 + y^3 + z^3);

i4 : isFPure R

o4 = true

i5 : S = ZZ/2[x,y,z];

i5 : isFPure ideal(y^2 - x^3)

o5 = false

i6 : isFPure ideal(z^2 - x*y*z + x*y^2 + x^2*y)

o6 = true
\end{MyVerbatim}
}\medskip

The algorithm works by applying Fedder's Criterion {\cite{FedderFPureRat}}, which states that a
local ring $(R, \mathfrak{m})$ is $F$-pure if and only if $(I^{[p]} : I)
\not\subseteq \mathfrak{m}^{[p]}$.
When \texttt{AtOrigin} is set to \texttt{true}, the algorithm checks $F$-purity only at the origin, by explicitly checking the above
containment. When \texttt{AtOrigin} is set to \texttt{false}, which is its default value, the
algorithm computes the non $F$-pure locus, by applying \texttt{frobeniusRoot}
to $I^{[p]} :I$. If the non $F$-pure locus is the whole ring, the algorithm
returns \texttt{true}.

\subsection{$F$-rationality}

\begin{definition}
Suppose that $R$ is a Cohen--Macaulay ring and that $T_{R} :
\omega_{R^{1/p}} \rightarrow \omega_{R}$ is the canonical dual to the
Frobenius. We say that $R$ has \emph{$F$-rational singularities}, of simply that \emph{$R$ is $F$-rational}, if there
are no non-zero proper submodules $M$ of $\omega_{R}$ such that
$T_{R}(M^{1/p}) \subseteq M$.
\end{definition}

The command \texttt{isFRational} checks if a ring is $F$-rational.

\medskip
{\small\setstretch{.67}
\begin{MyVerbatim}
i1 : S = ZZ/3[a,b,c,d,t];

i2 : M = matrix{{ a^2 + t^4, b, d }, { c, a^2, b^3 - d }};

             2       3
o2 : Matrix S  <--- S

i3 : I = minors(2, M);

o3 : Ideal of S

i4 : R = S/I;

i5 : isFRational R

o5 = true
\end{MyVerbatim}
}\medskip

The algorithm used by \texttt{isFRational} first checks if the ring is Cohen--Macaulay, unless the option \texttt{AssumeCM} is set to \texttt{true}. If the ring
is not Cohen--Macaulay then \texttt{false} is returned. Next, the algorithm
computes the test module $M \subseteq \omega_{R}$ and checks to see if
$\omega_{R} \subseteq M$;  see the next section for the description of a test module.

The options \texttt{AssumeDomain} and \texttt{FrobeniusRootStrategy} can be
used to improve the speed of the \texttt{testModule} computation. By default,
these options are set to \texttt{false} and \texttt{Substitution}, respectively.
Finally, if \texttt{AtOrigin} is set to \texttt{true}, then $F$-rationality is checked only at the origin.

\section{Test ideals}\label{Section: Test Ideals}

In this section, we explain how to compute parameter test modules, parameter test ideals, test ideals, and HSLG modules\footnote{HSLG modules can be used to give a scheme structure to the $F$-injective or $F$-pure locus.}.

\subsection{Parameter test modules}
%
%
Given a reduced ring $R$ of finite type over a perfect field $k$, the Frobenius map $R \hookrightarrow R^{1/p^e}$ is dual to $T : \omega_{R^{1/p^e}} \to \omega_R$.  As in \autoref{subsec.DualToFrobeniusOnQuotientRings}, we can represent the canonical module $\omega_R \subseteq R$ as an ideal, we can write $R = S/I$, and so we can find an element $u \in S^{1/p^e}$ representing the map $T:\omega_{R^{1/p^e}} \to \omega_R$; see \autoref{subsec.DualToFrobeniusOnQuotientRings} or \cite{KatzmanParameterTestIdealOfCMRings}.

\begin{definition}
   The \emph{parameter test submodule} is the smallest submodule $M \subseteq \omega_R$ (and hence ideal of $R$, since $M \subseteq \omega_R \subseteq R$) that agrees generically with $\omega_R$ and that satisfies
   %
   %
\[
T (M^{1/p^e}) \subseteq M
\]
for some $e > 0$ (or equivalently for all $e > 0$).
\end{definition}

Using \emph{Macaulay2}, we can compute this using the \texttt{testModule} command as follows.

\medskip
{\small\setstretch{.67}
\begin{MyVerbatim}
i1 : R = ZZ/5[x,y,z]/(x^4 + y^4 + z^4);

i2 : N = testModule R;

i3 : N#0

             2             2        2
o3 = ideal (z , y*z, x*z, y , x*y, x )

o3 : Ideal of R

i4 : N#1

o4 = ideal 1

o4 : Ideal of R
\end{MyVerbatim}
}\medskip

The output of \texttt{testModule} is a sequence, consisting of three items:
\begin{enumerate}[(1)]
\item The test module itself, given as an ideal of $R$.
\item The canonical module that contains the test module, given as an ideal of $R$.  (Note the representation of the canonical module as an ideal is not unique, it is only unique up to isomorphism, hence it is important to return this module as well).
\item The element $u$ described above (not displayed above, since it takes a lot of space).
\end{enumerate}

Note since this ring is Gorenstein, the canonical module is simply represented as the unit ideal. Here is another example, where the ring is not Gorenstein.

\medskip
{\small\setstretch{.67}
\begin{MyVerbatim}
i1 : R = ZZ/5[x,y,z]/(y*z, x*z, x*y);

i2 : N = testModule R;

i3 : N#0

             2   2   2
o3 = ideal (z , y , x )

o3 : Ideal of R

i4 : N#1

o4 = ideal (y + 2z, x - z)

o4 : Ideal of R
\end{MyVerbatim}
}\medskip

We briefly explain how this is computed:  First, we find a \emph{test element}.
\begin{remark}[Computation of test elements]
\label{rem.ComputationOfTestElements}
We recall that, roughly, an element of the Jacobian ideal that is not contained in any minimal prime is a test element \cite{HochsterFoundations}.  We compute test elements by computing random partial derivatives (and linear combinations thereof) until we find an element that does not vanish at all the minimal primes.  This method is much faster than computing the entire Jacobian ideal.  If it is known that the ring is a domain, setting \texttt{AssumeDomain => true} can speed this up further.
\end{remark}

After the test element $c$ has been identified, we pull back the ideal $\omega_R$ to an ideal $J \subseteq S$.  Next, we compute the following ascending sequence of ideals where $u$ represents $T:\omega_{R^{1/p^e}} \to \omega_R$ as above:
\begin{equation}
{\arraycolsep=1.4pt
\label{eq.AscendIdealExplanation}
\begin{array}{lll}
J_0 & \colonequals  cJ \\
J_1 &\colonequals   J_0 + (u J_0)^{[1/p]} & \\
J_2 &\colonequals   J_1 + (u J_1)^{[1/p]} & = J_0 + (u J_0)^{[1/p]} + (u^{1+p} J_0)^{[1/p^2]} \\
J_3 &\colonequals    J_2 + (uJ_2)^{[1/p]} & = J_0 + (u J_0)^{[1/p]} + (u^{1+p} J_0)^{[1/p^2]} + (u^{1+p+p^2} J_0)^{[1/p^3]} \\
\vdots
\end{array}
}
\end{equation}
As soon as this ascending sequence of ideals stabilizes, we are done.
In fact, because this strategy is used in several contexts, the user can call it directly for a chosen ideal $J$ and $u$ with the function \texttt{ascendIdeal} (this is done for test ideals below).

We can also compute parameter test modules of pairs $(\omega_R, f^{t})$ with $t \in \mathbb{Q}_{\geq 0}$.
This is done by modifying the element $u$ when the denominator of $t$ is not divisible by $p$.
When $t$ has $p$ in its denominator, we rely on the fact (see \cite{BlickleMustataSmithDiscretenessAndRationalityOfFThresholds,SchwedeTuckerTestIdealFiniteMaps}) that
\[
\begin{array}{rcl}
T(\tau(\omega_R, f^a)) & = & \tau(\omega_R, f^{a/p})\\
{\tt frobeniusRoot(1, u*I_1)} & {\tt =} & {\tt I_2}
\end{array}
\]
where the second line roughly explains how this is accomplished internally.
Here ${\tt I_1}$ is $\tau(\omega_R, f^a)$ pulled back to $S$ and, likewise, $I_2$ defines $\tau(\omega_R, f^{a/p})$ modulo the defining ideal of $R$.

\begin{remark}[Optimizations in \texttt{ascendIdeal} and other \texttt{testModule} computations]
Throughout the computations described above, we very frequently use the following fact:
\[
(f^p \cdot J)^{[1/p]} = f \cdot (J^{[1/p]}).
\]
To access this enhancement, one should try to pass functions like \texttt{ascendIdeal} and \texttt{frobeniusRoot} the elements and their exponents (see the documentation).
In particular, when computing the $p^e$-th Frobenius root of an ideal of the form $f^n \cdot J$, one should not multiply out $f^n$ and $J$, nor even raise $f$ to the $n^\mathrm{th}$ power, but rather call \texttt{frobeniusRoot(e, n, f, J)}.
\end{remark}

\subsection{Parameter test ideals}

The parameter test ideal is simply the annihilator of $\omega_R/\tau(\omega_R)$.  In other words, it is
\[
( \tau(\omega_R) : \omega_R ).
\]
This can also be described as
\[
\bigcap_{I} (I^* : I),
\]
where $I$ varies over ideals defined by a partial system of parameters, and $I^*$ denotes its tight closure.  This latter description is not computable, however.

\begin{example}\label{Example: parameter test ideal}
In this example we repeat the calculation in \cite[\S 9]{KatzmanParameterTestIdealOfCMRings}.

\medskip
{\small\setstretch{.67}
\begin{MyVerbatim}
i1 : R = ZZ/2[a..e];

i2 : E = matrix {{a, b, b, e}, {d, d, c, a}};

             2       4
o2 : Matrix R  <--- R

i3 : I = minors(2, E);

o3 : Ideal of R

i4 : S = R/I;

i5 : J = parameterTestIdeal S

o5 = ideal (c + d, b, a)

o5 : Ideal of S

i6 : J = substitute(J, R);

o6 : Ideal of R

i7 : mingens(J + I)

o7 = | c+d b a de |

             1       4
o7 : Matrix R  <--- R
\end{MyVerbatim}
}\medskip

\end{example}

\subsection{Test ideals}

For an $F$-finite reduced ring $R = S/I$, where $S$ is a regular ring, the (big) test ideal\footnote{The notion of test ideals was originally introduced in \cite{HochsterHunekeTC1} in the context of tight closure in finitely generated modules, whereas our notion of test ideals arises from
tight closure in possibly non finitely generated (``big'') modules. Confusingly, big test ideals are included in (small) test ideals.}
of $R$ is the smallest ideal $J$ in $R$, not contained in any minimal prime, such that for all $e > 0$ and all $\phi \in \Hom_R(R^{1/p^e}, R)$ we have
\[ \phi(J^{1/p^e}) \subseteq J.\]
In the case that $R$ is Gorenstein,
$\Hom_R(R^{1/p^e}, R)$ is a cyclic $R^{1/p^e}$-module generated by
$\Phi_e$, which corresponds with the map $T$ above based on the identification
$\omega_R \cong R$ \cite{BlickleSchwedeSurveyPMinusE}.  More generally, if $R$ is $\bQ$-Gorenstein with index not divisible by $p$, then for at least sufficiently divisible $e > 0$, such a generating $\Phi_e$ still exists.

If such a $\Phi_e$ exists, it can be identified to a generator of the module $(I^{[p^e]} : I) / I^{[p^e]}$ by Fedder's Lemma \cite{FedderFPureRat}.  Hence we can find a corresponding\footnote{This is done via the function \texttt{QGorensteinGenerator}. } $u \in I^{[p^e]} : I$.  In this case, if $c \in S$ is the pre-image of a test element of $R$, then setting $I_0 = cR$, it follows that
\[
\tau(R) = {\tt ascendIdeal(e, u, I_0)},
\]
where \texttt{ascendIdeal} is the method explained above, in \autoref{eq.AscendIdealExplanation}.
The $e$ here means all Frobenius roots are taken as multiples of $e$.  In other words, we first compute ${\tt I_0 +  (u \cdot I_0)^{1/p^{[e]}}}$.  Then we compute
\[
{\tt I_0 +  (u\cdot I_0)^{1/p^{[e]}} + (u^{p^e + 1} \cdot I_0)^{1/p^{[2e]}} }
\]
etc.

Here is an example (a $\bZ/3\bZ$-quotient, where $3 | (7-1)$), where exactly this logic occurs.

\medskip
{\small\setstretch{.67}
\begin{MyVerbatim}
i1 : T = ZZ/7[x,y];

i2 : S = ZZ/7[a,b,c,d];

i3 : f = map(T, S, {x^3, x^2*y, x*y^2, y^3});

i4 : I = ker f;

i5 : R = S/I;

i6 : testIdeal R

o6 = ideal 1

o6 : Ideal of R
\end{MyVerbatim}
}\medskip

However, the term $u$ constructed above can be quite complicated if $e > 1$ (which happens exactly when $(p -1)K_R$ is not Cartier).  For instance, even in the above example we have an extremely complex $u$:

\medskip
{\small\setstretch{.8}
\begin{MyVerbatim}
i7 : toString QGorensteinGenerator(1, R)

o7 = a^2*b^6*c^12+a^3*b^4*c^13+a^3*b^5*c^11*d+a^4*b^3*c^12*d+a^5*b*c^13*d+b^
     12*c^6*d^2+a^3*b^6*c^9*d^2+a^4*b^4*c^10*d^2+a^5*b^2*c^11*d^2+a^6*c^12*d
     ^2+b^13*c^4*d^3+a*b^11*c^5*d^3+a^2*b^9*c^6*d^3+a^4*b^5*c^8*d^3+a^5*b^3*
     c^9*d^3+a^6*b*c^10*d^3+a*b^12*c^3*d^4+a^2*b^10*c^4*d^4+a^3*b^8*c^5*d^4+
     a^4*b^6*c^6*d^4+a^5*b^4*c^7*d^4+a^6*b^2*c^8*d^4+a^7*c^9*d^4+a*b^13*c*d^
     5+a^2*b^11*c^2*d^5+a^3*b^9*c^3*d^5+a^4*b^7*c^4*d^5+a^5*b^5*c^5*d^5+a^6*
     b^3*c^6*d^5+a^7*b*c^7*d^5+a^2*b^12*d^6+a^3*b^10*c*d^6+a^4*b^8*c^2*d^6+a
     ^5*b^6*c^3*d^6+a^6*b^4*c^4*d^6+a^7*b^2*c^5*d^6+a^8*c^6*d^6+a^4*b^9*d^7+
     a^5*b^7*c*d^7+a^6*b^5*c^2*d^7+a^7*b^3*c^3*d^7+a^8*b*c^4*d^7+a^6*b^6*d^8
     +a^7*b^4*c*d^8+a^8*b^2*c^2*d^8+a^9*c^3*d^8+a^8*b^3*d^9+a^9*b*c*d^9+a^10
     *d^10
\end{MyVerbatim}
}\medskip

Therefore, we use a different strategy if either $(p-1)K_R$ is not Cartier or, more generally, if $R$ is $\bQ$-Gorenstein of index divisible by $p$.
In these situations, this alternate strategy typically appears to be faster.  We rely on the observation (see \cite{BlickleSchwedeTuckerTestAlterations}) that
\[
\tau(\omega_R, K_R) \cong \tau(R).
\]
In fact, by embedding $\omega_R \subseteq R$, we can compute $g$ so that $\tau(\omega_R, K_R) = g\cdot \tau(R)$.  We can therefore find $\tau(R)$ if we can find $\tau(\omega_R, K_R)$.
Next, if $K_R$ is $\bQ$-Cartier with $nK_R = \Div_R(f)$ for some $f \in R$ and $n > 0$, then
\[
\tau(\omega_R, K_R) =\tau(\omega_R, f^{1/n}).
\]
Thus we directly compute $\tau(\omega_R, f^{1/n})$ via the command \texttt{testModule(1/n, f)}.  Consider the following example, a $\mu_3$-quotient, which uses the logic described above.

\medskip
{\small
\setstretch{.67}
\begin{MyVerbatim}
i1 : T = ZZ/3[x,y];

i2 : S = ZZ/3[a,b,c,d];

i3 : f = map(T, S, {x^3, x^2*y, x*y^2, y^3});

o3 : RingMap T <--- S

i4 : I = ker f;

o4 : Ideal of S

i5 : R = S/I;

i6 : testIdeal R

o7 = ideal 1

o7 : Ideal of R
\end{MyVerbatim}
}\medskip

\begin{remark}[Non-graded caveats]
   It frequently happens that $(I^{[p^e]} : I)/I^{[p^e]}$ is principal but \emph{Macaulay2} cannot identify the principal generator (since \emph{Macaulay2} cannot always find minimal generators of non-graded ideals or modules).
   The same thing can happen when computing the element $u$ corresponding to the map $T : \omega_{R^{1/p}} \to \omega_R$, as described in \autoref{subsec.DualToFrobeniusOnQuotientRings}.
   In such situations, instead of a single $u$, \emph{Macaulay2} will produce $u_1, \dots, u_n$ (all multiples of $u$, and $u$ is a linear combination of the $u_i$).  Instead of computing the ideal
\[
(u \cdot J)^{[1/p]},
\]
we compute
\[
(u_1 \cdot J)^{[1/p]} + \dots + (u_n \cdot J)^{[1/p]},
\]
which will produce the same answer.
\end{remark}

We can similarly use the \texttt{testIdeal} command to compute test ideals of pairs, $\tau(R, f^t)$, and even mixed test ideals, $\tau(R, f_1^{t_1} \cdots f_n^{t_n})$.

\subsection{HSLG module; computing $F$-pure submodules of rank-1 Cartier modules}

Again, consider the maps $T^e : \omega_{R^{1/p^e}} \to \omega_R$ we have considered throughout this section.  It is a theorem of Hartshorne-Speiser, Lyubeznik, and Gabber \cite{HartshorneSpeiserLocalCohomologyInCharacteristicP,LyubeznikFModulesApplicationsToLocalCohomology,Gabber.tStruc} that the descending images
\begin{equation}\label{eq: descending images}
\omega_R \supseteq T(\omega_{R^{1/p}}) \supseteq \dots \supseteq T^e(\omega_{R^{1/p^e}}) \supseteq T^{e+1}(\omega_{R^{1/p^{e+1}}}) \supseteq \cdots
\end{equation}
stabilize for $e \gg 0$.
The function \texttt{FPureModule} computes the stable submodule in this chain, called the \emph{HSLG module}, and returns a sequence containing the following items:
\begin{enumerate}[(1)]
   \item The HSLG module;
   \item The canonical module, embedded (non-uniquely) as an ideal;
   \item The element $u$ representing the map on the canonical module (see
   \autoref{subsec.DualToFrobeniusOnQuotientRings});
   \item The value of $e > 0$ at which the descending images in \eqref{eq: descending images} stabilize.  This is sometimes also called the HSLG number of the canonical module as a Cartier module.
\end{enumerate}

\medskip
{\small
\setstretch{.67}
\begin{MyVerbatim}
i1 : R = ZZ/3[x,y,z]/(x^3 + y^4 + z^5);

i2 : L = FPureModule R;

i3 : L#0

                       2        2   3
o3 = ideal (y*z, x*z, y , x*y, x , z )

o3 : Ideal of R

i4 : L#1

o4 = ideal 1

o4 : Ideal of R

i5 : L#3

o5 = 1

i6 : R = ZZ/3[x,y,z]/(x^3 + y^4 + z^5);

i7 : L = FPureModule R;

i8 : L#0

                       2        2   3
o8 = ideal (y*z, x*z, y , x*y, x , z )

o8 : Ideal of R

i9 : L#1

o9 = ideal 1

o9 : Ideal of R

i10 : L#3

o10 = 1
\end{MyVerbatim}
}\medskip


More generally, for any ideal $J$ with a map $\phi : J^{1/p^e} \to J$, we have that the images
\[
J \supseteq \phi(J^{1/p^e}) \supseteq \phi^2(J^{1/p^{2e}}) \supseteq \cdots
\]
stabilize as well (in fact, the analogous result even holds for modules \cite{Gabber.tStruc}).

\subsection{Computing the level of a polynomial}
Another interesting invariant that can be calculated using this package
is the so--called \emph{level} of a polynomial; more precisely:

\begin{definition}
Let $K$ be an $F$--finite field of prime characteristic $p,$ and let
$f\in R=K[x_1,\ldots, x_d].$ We define the \emph{level of} $f$ as the
smallest possible integer $e$ where the descending chain
\[
R=(f^{p^0-1})^{[1/p^0]}\supseteq (f^{p-1})^{[1/p]}\supseteq
(f^{p^2-1})^{[1/p^2]}\supseteq\ldots
\supseteq (f^{p^i-1})^{[1/p^i]}\supseteq\ldots
\]
stabilizes.
\end{definition}
This invariant was introduced in \cite{AlvarezBlickleLyubeznik2005}.  It is also essentially the same data as the HSLG number of the pair $(R, f^1)$ as computed in (4) above (it is that number $+1$).
One interesting particular
case is when $K=\mathbb{F}_p$ and $f$ is the defining homogeneous polynomial of a hyperelliptic curve
of genus $g;$ when $g=1,$ it was proved in
\cite{BoixDeStefaniVanzo2015} that the corresponding elliptic curve defined by $f$
is ordinary if and only if its level is $1,$ (equivalently, if and only if
$(R/(f))_{(x,y,z)}$ is $F$--pure, 
) and supersingular if and only if its level is $2.$ When
the genus is at least $2,$ level $2$ is a necessary (but not sufficient) condition for the curve
for being ordinary \cite{BlancoBoixFordhamYilmaz2018}; in this case, one also
has that, if the curve is supersingular, then its level has to be at least $3,$ so the level can
always distinguish these two properties in any genus. We illustrate these results by means
of the following examples.

\medskip
{\small
\setstretch{.67}
\begin{MyVerbatim}
i1 : R = ZZ/2[x,y,z];

i2 : f = x^3 + y^2*z + y*z^2;

i3 : frobeniusPower(1/2, ideal f)

o3 = ideal (z, y, x)

o3 : Ideal of R

i4 : u = f^3;

i5 : frobeniusPower(1/4, ideal u)

o5 = ideal (z, y, x)

o5 : Ideal of R

i6 : ((FPureModule(1, f))#3) + 1

o6 = 2
\end{MyVerbatim}
}\medskip

This shows that this elliptic curve has level $2,$ hence it is supersingular. The next example
shows a non-ordinary hyperelliptic curve of genus $2$ with level $2$.

\medskip
{\small
\setstretch{.67}
\begin{MyVerbatim}
i1 : R = ZZ/11[x,y,z];

i2 : f = y^2*z^3 - x^5 - 2*z^5;

i3 : frobeniusPower(1/11, ideal f^10)

             2        3
o3 = ideal (z , x*z, x )

o3 : Ideal of R

i4 : frobeniusPower(1/121, ideal f^120)

             2        3
o4 = ideal (z , x*z, x )

o4 : Ideal of R

i5 : ((FPureModule(1, f))#3) + 1

o5 = 2
\end{MyVerbatim}
}\medskip

\section{Ideals compatible with a given $p^{-e}$-linear map}\label{Section: compatible ideals}

Throughout this section, let $R$ denote a polynomial $\mathbb{K}[x_1, \dots, x_n]$.
In this section we address the following question:
given a $R$-linear map $\phi: R^{1/p^e} \rightarrow R$, what are the ideals $I\subseteq R$ such that $\phi(I^{1/p^e})\subseteq I$?
We call these ideals $\phi$-compatible.

Recall that  $\Hom_R( R^{1/p^e}, R)$  is a principal $R^{1/p^e}$-module generated
the \emph{trace map} $T\in \Hom_R(R^{1/p^e}, R)$, constructed as in \autoref{ss: math behind}, (cf. \cite[Lemma 1.6]{FedderFPureRat} and \cite[Example 1.3.1]{BrionKumarFrobeniusSplitting}).
%
%

We can now write our given $\phi$ as multiplication by some $u^{1/p^e}$ followed by $T$ and it is not hard to see that
an ideal $I\subseteq R$ is $\phi$-compatible if and only if $u I \subseteq I^{[p^e]}$.

\begin{theorem}\label{Theorem: finitely many compatible primes}
If $\phi$ is surjective, there are finitely many $\phi$-compatible ideals, consisting of all possible intersections
of $\phi$-compatible prime ideals \textup(cf.\ \cite{KumarMehtaFiniteness}, \cite{SchwedeFAdjunction},
\cite{SharpGradedAnnihilatorsOfModulesOverTheFrobeniusSkewPolynomialRing}, \cite{EnescuHochsterTheFrobeniusStructureOfLocalCohomology}\textup).
In general, there are finitely many $\phi$-compatible prime ideals not containing $\sqrt{\Image \phi}$ \textup(cf.\ \cite{KatzmanSchwedeAlgorithm}\textup).

\end{theorem}

The method \texttt{compatibleIdeals} produces the finite set of $\phi$-compatible prime ideals in the second claim of \autoref{Theorem: finitely many compatible primes}.

\medskip
{\small
\setstretch{.67}
\begin{MyVerbatim}
i1 : R = ZZ/3[u,v];

i2 : u = u^2*v^2;

i3 : compatibleIdeals u

o3 = {ideal v, ideal (v, u), ideal u}

o3 : List
\end{MyVerbatim}
}\medskip

The defining condition $u I \subseteq I^{[p^e]}$ for $\phi$-compatible ideals allows us to
think of these in a dual form: write $\mathfrak{m}=(x_1, \dots, x_n)$,
$E=E_{R_\mathfrak{m}}(R_{\mathfrak{m}}/\mathfrak{m})=H^n_{\mathfrak{m}} (R)$, and let
$\Theta: E \rightarrow E$ be the $p^e$-linear map\footnote{That is, $\Theta$ is additivive and $\Theta (r a)= r^{p^e} \Theta (a)$ for all $a\in E$ and $r\in R$.}
induced from the Frobenius map on $R$.
If $\psi=u \Theta$, then $\psi \Ann_E I \subseteq \Ann_E I$ if and only if $u I \subseteq  I^{[p^e]}$.
Thus finding all $R$-submodules of $E$ compatible with $\psi=u \Theta$ also amounts to finding all
$\phi$-compatible ideals, where $\phi=T \circ u^{1/p^e}$.

\begin{example}
We return to \autoref{Example: parameter test ideal}.
In (cf.~\cite[\S 9]{KatzmanParameterTestIdealOfCMRings}) it is shown that
there is a surjection $\Ann_E I \rightarrow H^2_{\mathfrak{m}} (R/I)$
which is compatible  with the induced $p^1$-linear map on $H^2_{\mathfrak{m}} (R/I)$
and the $p^1$-linear map $u \Theta$ on $\Ann_E I$, where $u$ is computed as follows.

\medskip
{\small
\setstretch{.67}
\begin{MyVerbatim}
i1 : R = ZZ/2[a..f];

i2 : E = matrix {{a, b, b, e}, {d, d, c, a}};

             2       4
o2 : Matrix R  <--- R

i3 : I = minors(2, E);

o3 : Ideal of R

i4 : S = R/I;

i5 : isCohenMacaulay S

o5 = true
\end{MyVerbatim}
}\medskip

In \cite{KatzmanParameterTestIdealOfCMRings} it is shown that as $R/I$ is Cohen--Macaulay, $u$
can be taken as the generator of the cyclic module $(I^{[p]}:I) \cap (\Omega^{[p]}:\Omega)/I^{[p]}$
where $\Omega \subseteq R$ is an ideal whose image in $R/I$ is a canonical module of that ring.

\medskip
{\small
\setstretch{.67}
\begin{MyVerbatim}
i6 : omega = canonicalIdeal S

o6 = ideal (e, d, a)

o6 : Ideal of S

i7 : omega = substitute(omega, R) + I;

o7 : Ideal of R

i8 : u = intersect((frobenius I):I, (frobenius omega) : omega);

o8 : Ideal of R

i9 : u = compress((gens u) % gens(frobenius I))

o9 = | a3bc+a3bd+a2cde+abcde+abd2e+b2d2e+cd2e2+d3e2 |

             1       1
o9 : Matrix R  <--- R

i10 :  u = first first entries u

       3       3       2                         2     2 2       2 2    3 2
o10 = a b*c + a b*d + a c*d*e + a*b*c*d*e + a*b*d e + b d e + c*d e  + d e

o10 : R
\end{MyVerbatim}
}
\medskip

Now we can compute all annihilators of $R$-submodules of $E$ stable under the $p^1$-linear map $uT$.

\medskip
{\small
\setstretch{.67}
\begin{MyVerbatim}
i11 : L = compatibleIdeals u;

i12 : print \ L;

       3       3       2                         2     2 2       2 2    3 2
ideal(a b*c + a b*d + a c*d*e + a*b*c*d*e + a*b*d e + b d e + c*d e  + d e )

               2
ideal (a + b, a  + d*e)

ideal (e, a, d)

ideal (e, d, a, c)

ideal (e, d, c, b, a)

ideal (e, d, b, a)

ideal (e, b, a)

ideal (e, b, a, c + d)

                      2
ideal (c + d, a + b, b  + d*e)

ideal (d, a)

ideal (d, a, c)

ideal (d, c, b, a)

ideal (d, b, a)
\end{MyVerbatim}
}
\medskip

We can also compute  all annihilators of $R$-submodules of $H^2_{\mathfrak{m}} (R/I)$ stable under the $p^1$-linear map induced from the Frobenius map on $R/I$.

\medskip
{\small
\setstretch{.67}
\begin{MyVerbatim}
i13 : print \ unique apply(L, J -> J : omega);

       3       3       2                         2     2 2       2 2    3 2
ideal(a b*c + a b*d + a c*d*e + a*b*c*d*e + a*b*d e + b d e + c*d e  + d e )

               2
ideal (a + b, b  + d*e)

ideal (e, d, a)

ideal 1

ideal (e, b, a)

ideal (e, c + d, b, a)

                      2
ideal (c + d, a + b, b  + d*e)

ideal (d, a)

ideal (d, c, a)

ideal (d, c, b, a)

ideal (d, b, a)
\end{MyVerbatim}
}
\medskip

\end{example}

\section{Future plans}

In \cite{KatzmanZhangAlgorithm}, the algorithms behind the methods in \autoref{Section: compatible ideals} were extended
to compute prime annihilators of submodules of Artinian modules compatible with a given $p^{e}$-linear map.
This would require, among other things, a faster implementation of our method for finding Frobenius roots of submodules of free modules.

On the other hand, it should be possible to compute test ideals of pairs $(R, \fra^t)$ even when $\fra$ is not principal.  One strategy to do this is outlined in \cite{SchwedeTuckerTestIdealsOfNonPrincipal} although certain improvement can be made.

We hope to achieve all these things during a future \emph{Macaulay2} workshop.  We also want to bring to the reader's attention the package \emph{FThresholds}, which computes $F$-pure thresholds, $F$-thresholds, $F$-jumping numbers and more!

\bibliographystyle{skalpha}
\bibliography{MainBib}

\def\cfudot#1{\ifmmode\setbox7\hbox{$\accent"5E#1$}\else
  \setbox7\hbox{\accent"5E#1}\penalty 10000\relax\fi\raise 1\ht7
  \hbox{\raise.1ex\hbox to 1\wd7{\hss.\hss}}\penalty 10000 \hskip-1\wd7\penalty
  10000\box7}
\providecommand{\bysame}{\leavevmode\hbox to3em{\hrulefill}\thinspace}
\providecommand{\MR}{\relax\ifhmode\unskip\space\fi MR}
\providecommand{\MRhref}[2]{%
  \href{http://www.ams.org/mathscinet-getitem?mr=#1}{#2}
}
\providecommand{\href}[2]{#2}
\begin{thebibliography}{BCBFY18}

\bibitem[{\`A}MBL05]{AlvarezBlickleLyubeznik2005}
{\sc J.~{\`A}lvarez~Montaner, M.~Blickle, and G.~Lyubeznik}: \emph{Generators
  of {$D$}-modules in positive characteristic}, Math. Res. Lett. \textbf{12}
  (2005), no.~4, 459--473. {\sf\scriptsize 2155224}

\bibitem[BCBFY18]{BlancoBoixFordhamYilmaz2018}
{\sc I.~Blanco-Chac\'on, A.~F. Boix, S.~Fordham, and E.~S. Yilmaz}:
  \emph{Differential operators and hyperelliptic curves over finite fields},
  Finite Fields Appl. \textbf{51} (2018), 351--370. {\sf\scriptsize 3781411}

\bibitem[BMS08]{BlickleMustataSmithDiscretenessAndRationalityOfFThresholds}
{\sc M.~Blickle, M.~Musta{\c{t}}{\v{a}}, and K.~E. Smith}: \emph{Discreteness
  and rationality of {$F$}-thresholds}, Michigan Math. J. \textbf{57} (2008),
  43--61, Special volume in honor of Melvin Hochster.

\bibitem[BMS09]{BlickleMustataSmithFThresholdsOfHypersurfaces}
{\sc M.~Blickle, M.~Musta{\c{t}}{\u{a}}, and K.~E. Smith}:
  \emph{{$F$}-thresholds of hypersurfaces}, Trans. Amer. Math. Soc.
  \textbf{361} (2009), no.~12, 6549--6565.

\bibitem[BS13]{BlickleSchwedeSurveyPMinusE}
{\sc M.~Blickle and K.~Schwede}: \emph{$p^{-1}$-linear maps in algebra and
  geometry}, Commutative Algebra, Springer, 2013, pp.~123--205.

\bibitem[BSTZ10]{BlickleSchwedeTakagiZhang}
{\sc M.~Blickle, K.~Schwede, S.~Takagi, and W.~Zhang}: \emph{Discreteness and
  rationality of {$F$}-jumping numbers on singular varieties}, Math. Ann.
  \textbf{347} (2010), no.~4, 917--949.

\bibitem[BST15]{BlickleSchwedeTuckerTestAlterations}
{\sc M.~Blickle, K.~Schwede, and K.~Tucker}: \emph{{$F$}-singularities via
  alterations}, Amer.\ J.\ Math. \textbf{137} (2015), no.~1, 61--109.

\bibitem[BDSV15]{BoixDeStefaniVanzo2015}
{\sc A.~F. Boix, A.~De~Stefani, and D.~Vanzo}: \emph{An algorithm for
  constructing certain differential operators in positive characteristic},
  Matematiche (Catania) \textbf{70} (2015), no.~1, 239--271. {\sf\scriptsize
  3351468}

\bibitem[BK05]{BrionKumarFrobeniusSplitting}
{\sc M.~Brion and S.~Kumar}: \emph{Frobenius splitting methods in geometry and
  representation theory}, Progress in Mathematics, vol. 231, Birkh\"auser
  Boston Inc., Boston, MA, 2005.

\bibitem[EH08]{EnescuHochsterTheFrobeniusStructureOfLocalCohomology}
{\sc F.~Enescu and M.~Hochster}: \emph{The {F}robenius structure of local
  cohomology}, Algebra Number Theory \textbf{2} (2008), no.~7, 721--754.

\bibitem[Fed83]{FedderFPureRat}
{\sc R.~Fedder}: \emph{{$F$}-purity and rational singularity}, Trans. Amer.
  Math. Soc. \textbf{278} (1983), no.~2, 461--480.

\bibitem[Gab04]{Gabber.tStruc}
{\sc O.~Gabber}: \emph{Notes on some {$t$}-structures}, Geometric aspects of
  Dwork theory. Vol. I, II, Walter de Gruyter GmbH \& Co. KG, Berlin, 2004,
  pp.~711--734.

\bibitem[GS]{M2}
{\sc D.~R. Grayson and M.~E. Stillman}: \emph{Macaulay2, a software system for
  research in algebraic geometry}.

\bibitem[HW02]{HaraWatanabeFRegFPure}
{\sc N.~Hara and K.-I. Watanabe}: \emph{{$F$}-regular and {$F$}-pure rings vs.\
  log terminal and log canonical singularities}, J. Algebraic Geom. \textbf{11}
  (2002), no.~2, 363--392.

\bibitem[HS77]{HartshorneSpeiserLocalCohomologyInCharacteristicP}
{\sc R.~Hartshorne and R.~Speiser}: \emph{Local cohomological dimension in
  characteristic {$p$}}, Ann. of Math. (2) \textbf{105} (1977), no.~1, 45--79.

\bibitem[Her14]{HernandezFPureThresholdOfBinomial}
{\sc D.~J. Hern\'andez}: \emph{{$F$}-pure thresholds of binomial
  hypersurfaces}, Proc. Amer. Math. Soc. \textbf{142} (2014), no.~7,
  2227--2242.

\bibitem[Her15]{HernandezFInvariantsOfDiagonalHyp}
{\sc D.~J. Hern\'andez}: \emph{{$F$}-invariants of diagonal hypersurfaces},
  Proc.\ Amer.\ Math.\ Soc. \textbf{143} (2015), no.~1, 87--104.

\bibitem[HT17]{HernandezTeixeiraFThresholdFunctions}
{\sc D.~J. Hern\'andez and P.~Teixeira}: \emph{{$F$}-threshold
  functions\textup: {S}yzygy gap fractals and the two-variable homogeneous
  case}, J.\ Symbolic Comput. \textbf{80} (2017), 451--483.

\bibitem[HTW18a]{HernandezTeixeiraWittFrobeniusPowers}
{\sc D.~J. Hern\'andez, P.~Teixeira, and E.~E. Witt}: \emph{Frobenius powers},
  preprint, \href{https://arxiv.org/abs/1802.02705}{arXiv:1802.02705
  [math.AC]}, 2018.

\bibitem[HTW18b]{hernandez+etal.frobenius_examples}
{\sc D.~J. Hern\'andez, P.~Teixeira, and E.~E. Witt}: \emph{Frobenius powers of
  some monomial ideals}, to appear in J. Pure Appl. Algebra,
  \href{https://arxiv.org/abs/1808.09508}{arXiv:1808.09508 [math.AC]}, 2018.

\bibitem[Hoc07]{HochsterFoundations}
{\sc M.~Hochster}: \emph{Foundations of tight closure theory}, lecture notes
  from a course taught at the University of Michigan Fall 2007 (2007).

\bibitem[HH90]{HochsterHunekeTC1}
{\sc M.~Hochster and C.~Huneke}: \emph{Tight closure, invariant theory, and the
  {B}rian\c con-{S}koda theorem}, J. Amer. Math. Soc. \textbf{3} (1990), no.~1,
  31--116.

\bibitem[Kat08]{KatzmanParameterTestIdealOfCMRings}
{\sc M.~Katzman}: \emph{Parameter-test-ideals of {C}ohen-{M}acaulay rings},
  Compos. Math. \textbf{144} (2008), no.~4, 933--948.

\bibitem[Kat10]{KatzmanFrobeniusMapsOnInjectiveHulls}
{\sc M.~Katzman}: \emph{Frobenius maps on injective hulls and their
  applications to tight closure}, J. Lond. Math. Soc. \textbf{81} (2010),
  no.~3, 589--607.

\bibitem[KLZ09]{KatzmanLyubeznikZhangOnDiscretenessAndRationality}
{\sc M.~Katzman, G.~Lyubeznik, and W.~Zhang}: \emph{On the discreteness and
  rationality of {$F$}-jumping coefficients}, J. Algebra \textbf{322} (2009),
  no.~9, 3238--3247.

\bibitem[KS12]{KatzmanSchwedeAlgorithm}
{\sc M.~Katzman and K.~Schwede}: \emph{An algorithm for computing compatibly
  {F}robenius split subvarieties}, J. Symbolic Comput. \textbf{47} (2012),
  no.~8, 996--1008.

\bibitem[KZ14]{KatzmanZhangAlgorithm}
{\sc M.~Katzman and W.~Zhang}: \emph{Annihilators of artinian modules
  compatible with a {F}robenius map}, J. Symbolic Comput. \textbf{60} (2014),
  29--46.

\bibitem[KM09]{KumarMehtaFiniteness}
{\sc S.~Kumar and V.~B. Mehta}: \emph{Finiteness of the number of compatibly
  split subvarieties}, Int. Math. Res. Not. IMRN (2009), no.~19, 3595--3597.

\bibitem[Kun69]{KunzCharacterizationsOfRegularLocalRings}
{\sc E.~Kunz}: \emph{Characterizations of regular local rings for
  characteristic {$p$}}, Amer. J. Math. \textbf{91} (1969), 772--784.

\bibitem[Lyu97]{LyubeznikFModulesApplicationsToLocalCohomology}
{\sc G.~Lyubeznik}: \emph{{$F$}-modules: applications to local cohomology and
  {$D$}-modules in characteristic {$p>0$}}, J. Reine Angew. Math. \textbf{491}
  (1997), 65--130.

\bibitem[Rai]{PushForward}
{\sc C.~Raicu}: \emph{{P}ush{F}orward.m2}.

\bibitem[Sch09]{SchwedeFAdjunction}
{\sc K.~Schwede}: \emph{{$F$}-adjunction}, Algebra Number Theory \textbf{3}
  (2009), no.~8, 907--950.

\bibitem[ST12]{SchwedeTuckerTestIdealSurvey}
{\sc K.~Schwede and K.~Tucker}: \emph{A survey of test ideals}, Progress in
  Commutative Algebra 2. Closures, Finiteness and Factorization (C.~Francisco,
  L.~C. Klinger, S.~M. Sather-Wagstaff, and J.~C. Vassilev, eds.), Walter de
  Gruyter GmbH \& Co. KG, Berlin, 2012, pp.~39--99.

\bibitem[ST14a]{SchwedeTuckerTestIdealFiniteMaps}
{\sc K.~Schwede and K.~Tucker}: \emph{On the behavior of test ideals under
  finite morphisms}, J.\ Algebraic Geom. \textbf{23} (2014), no.~3, 399--443.

\bibitem[ST14b]{SchwedeTuckerTestIdealsOfNonPrincipal}
{\sc K.~Schwede and K.~Tucker}: \emph{Test ideals of non-principal ideals:
  Computations, jumping numbers, alterations and division theorems}, J. Math.
  Pures Appl. \textbf{102} (2014), no.~5, 891 -- 929.

\bibitem[Sha07]{SharpGradedAnnihilatorsOfModulesOverTheFrobeniusSkewPolynomialRing}
{\sc R.~Y. Sharp}: \emph{Graded annihilators of modules over the {F}robenius
  skew polynomial ring, and tight closure}, Trans. Amer. Math. Soc.
  \textbf{359} (2007), no.~9, 4237--4258 (electronic).

\end{thebibliography}

\end{document}